\documentclass[twoside]{amsart}
\usepackage{amsmath,amssymb,amsthm}

\newcommand{\N}{\mathbf{N}}
\newcommand{\Z}{\mathbf{Z}}

\newcommand{\R}{\mathbf{R}}

\newcommand{\SL}{\mathrm{SL}}
\newcommand{\PSL}{\mathrm{PSL}}

\newcommand{\lcm}{\mathrm{lcm}}

\newcommand{\Id}{\mathrm{Id}}

\newcommand{\T}{\mathrm{T}}

\theoremstyle{definition}

\title{Veech surfaces associated with rational billiards}
\author{Samuel Leli\`evre} 
\date{May 23, 2002}

\begin{document}

\begin{abstract}
A nice trick for studying the billiard flow in a rational polygon 
is to unfold the polygon along the trajectories. This gives rise to
a translation or half-translation surface tiled by the original polygon,
or equivalently an Abelian or quadratic differential.

Veech surfaces are a special class of translation surfaces with a 
large group of affine automorphisms, and interesting dynamical 
properties. 
The first examples of Veech surfaces came from rational billiards.

We first present the mathematical objects and fix some vocabulary and 
notation. Then we review known results about Veech surfaces arising 
from rational billiards. The interested reader will find 
annex tables on the author's web page.
\end{abstract}

\maketitle

\section{The objects}
\subsection{Polygons}
We consider billiards in Euclidean (plane) polygons.
A polygon is called rational if the angles made by its sides are all 
rational in unit $\pi$.
(Note concerning non simply connected polygons: this condition is not 
only for adjacent sides.) A rational billiard is a billiard in 
a rational polygon.
All billiards in this paper are rational.

A triangle is determined by its angles up to a scaling 
factor. Thus the parameter space for rational triangles is discrete.
This is not the case for other polygons.

We denote by $\T(a,b,c)$ the triangle 
$(\cfrac{a\,\pi}{d},\cfrac{b\,\pi}{d},\cfrac{c\,\pi}{d})$ where
$a,b,c\in \N^*$ and $d=a+b+c$. 
This notation is unique if we require $\gcd(a,b,c)=1$ and 
$a\leqslant b \leqslant c$.

A triangle is acute if all its angles are acute, right if it has a 
right angle, obtuse if it has an obtuse angle.

A triangle is equilateral if its three angles are equal, isosceles if 
two of its angles are equal, scalene otherwise.

\subsection{Translation and half-translation surfaces}
A translation structure on a surface is an atlas on this surface with 
finitely many punctures, whose transition functions are translations. 
A surface which admits a translation structure is called a translation 
surface.

A half-translation structure on a surface is an atlas on this surface 
with finitely many punctures, whose transition functions are 
translations possibly composed with a central symmetry. A surface which 
admits a half-translation structure is called a half-translation 
surface.

A translation surface is a particular case of a half-translation 
surface. 
A half-translation surface is not a translation surface in general. 
But from a half-translation surface one can always obtain a translation 
surface by a branched double cover.

There is a correspondance between translation surfaces and Abelian 
differentials, and between half-translation surfaces and quadratic 
differentials.

\subsection{Surfaces associated with a rational billiard}

\subsubsection*{Translation surface}
There is a classical construction of a translation surface from a 
rational billiard (see \cite{FoxKsh36, KtkZml75}). 
Starting with a polygon with angles 
$(m_i/n_i)\pi$, and letting $N=\lcm(n_i)$, 
to each vertex of the 
polygon correspond $N/n_i$ singular points on the surface with 
multiplicity $m_i$ (cone points with angle $m_i2\pi$ or, in other words, 
zeros of order $m_i-1$ of the Abelian differential).

\subsubsection*{Half-translation surface}
One can consider a reduced construction in order to obtain a 
half-translation surface. Starting with a polygon with angles 
$(m_i/n_i)(\pi/2)$, and letting $N=\lcm(n_i)$, 
to each vertex of the polygon correspond $N/n_i$ singular points 
on the surface with multiplicity $m_i/2$ (cone points with angle 
$m_i\pi$ or, in other words, zeros of order $m_i-2$ of the 
quadratic differential).

\subsubsection*{Comparison}
Both constructions are the same if the classical $N$ is odd. Otherwise 
the translation surface is a branched double cover of the 
half-translation surface.

\subsubsection*{Remark}
For each $k$ we can construct a 
``$1/k$-translation'' surface, endowed with a differential of order $k$. 
In particular, it is a flat surface with cone type singularities, 
whose angles are multiples of $2\pi/k$. 
The information for cone points can be obtained from angles of the 
polygon as above, by expressing them in unit $\pi/k$.

In particular, the minimal construction consists in gluing just two copies 
of the initial polygon.

However, translation and half-translation surfaces (or Abelian and 
quadratic differentials) play a distinguished role because there is an 
action of $\SL(2,\R)$ on the space of these surfaces (or of these 
differentials).

\subsection{Veech group}
On a translation surface we can consider affine diffeomorphisms 
(diffeomorphisms on the surface punctured at its singularities, that 
are affine in the charts of the translation structure, and extend to 
a homeomorphism of the whole surface). The derivatives (or linear 
parts) of these diffeomorphisms form a subgroup of $\SL^\pm(2,\R)$. 
If we consider only those diffeomorphisms which preserve orientation, 
we get a subgroup of $\SL(2,\R)$.

For half-translation surfaces we work modulo $-\Id$ so that we obtain 
subgroups of $\PSL^\pm(2,\R)$ or $\PSL(2,\R)$.

We call Veech group the subgroup of $\PSL(2,\R)$ obtained for the 
translation surface with only non-removable singularities (those that 
correspond to zeros of the Abelian differential).

\subsubsection*{Triangle groups}
We denote by $\triangle(p,q,r)$ the $(p,q,r)$-triangle group. See 
\cite{Brd83} for more information about these groups. This notation 
really refers to a conjugation class of subgroups of $\PSL(2,\R)$.
Note that the class of $\PSL(2,\Z)$ is $\triangle(2,3,\infty)$.

\subsubsection*{Lattice property}
A surface has the lattice property (or Veech property) if its Veech 
group has finite covolume in $\PSL(2,\R)$.
A polygon is said to have the lattice property if its associated 
translation surface does.

The lattice property is interesting in that it implies good ergodic 
properties, namely the Veech dichotomy: in each direction the 
flow is either periodic or uniquely ergodic.

\section{The results}

In this section we review polygons (mainly triangles) known to have 
or not to have the lattice property. We begin with the arithmetic 
examples, those that are the closest to the basic example of the 
square. We then move to regular polygons, and then to right and acute 
triangles for which the Veech property is characterized. 
For scalene triangles, only partial answers are known.

\subsection{Arithmetic case}

\begin{quote}
\begin{itshape}
The translation surface associated with the billiard in a rectangle, 
in a right isosceles triangle, in an equilateral triangle or in 
the $(\pi/6,\pi/3,\pi/2)$ triangle is a flat torus. Its Veech group
is $\triangle(2,3,\infty)$.
\end{itshape}
\end{quote}
In the following subsections we always skip the torus case.

Gutkin and Judge proved the following in \cite{GutJdg96,GutJdg00}.
\begin{quote}
\begin{itshape}
A translation surface has an arithmetic (conjugate to a subgroup of 
$\PSL(2,\Z)$) Veech group if and only if it 
is (translation) tiled by a Euclidean parallelogram, or in other 
words if it is  a (translation) cover of a one-punctured flat torus.
\end{itshape}
\end{quote}

\subsection{Regular polygons}
All regular polygons have the lattice property.

The Veech group of the translation surface associated with the regular 
$n$-gon for $n\geqslant 5$ is a subgroup of $\triangle(2,n,\infty)$. 
Its index can be given as follows \cite{Vch92}. Let $\varepsilon(n) = 
\gcd(2,n)$, let $N=n/\varepsilon(n)$, let $\sigma(n) = \gcd(4,n)$.

Let $\omega(n) = n 
\displaystyle{\prod_{\substack{p|n \\ p \text{ prime}}}(1 + 1/p)}$.
Then the index is $\cfrac{\omega(N)}{\omega(\sigma(N))}\,\varepsilon(n)$.

\subsection{Right triangles}

\begin{quote}
\begin{itshape}
  A right triangle has the lattice property if and only if
  its smallest angle is $\pi/n$ for some $n\geqslant4$.

  For $n\geqslant5$, the Veech group of the corresponding translation 
  surface is $\triangle(2,n,\infty)$ if $n$ is odd, and
  $\triangle(m,\infty,\infty)$ if $n=2m$. 
\end{itshape}
\end{quote}

Vorobets showed that the condition is sufficient and gave the
corresponding Veech group \cite[section 4]{Vrb96}.  Kenyon and Smillie
showed that the condition is necessary \cite[section 6]{KenSmi00}.

\subsection{Acute triangles}

\subsubsection{Isosceles}

\begin{quote}
\begin{itshape}
  An acute isosceles triangle has the lattice property if and only if 
  its apex angle is $\pi/n$ for some $n\geqslant3$.

  For $n\geqslant4$, the Veech group of the associated translation surface 
   is $\triangle(n,\infty,\infty)$.
\end{itshape}
\end{quote}

Showing that the condition is sufficient reduces to the right 
triangle case by an unfolding construction, see \cite[section 5]{Vrb96} 
or \cite{GutJdg00}.
Kenyon and Smillie showed that the condition is necessary 
\cite[section 6]{KenSmi00}.
The Veech group was given for an example in \cite{EarGrd97}, 
and for the general case in \cite{HbtSch00}.

\subsubsection{Scalene}

\begin{quote}
\begin{itshape}
An acute scalene triangle has the lattice property if and only if it is one
of the exceptional triangles
$$(\frac{2\pi}{9},\frac{\pi}{3},\frac{4\pi}{9}),\quad
(\frac{\pi}{4},\frac{\pi}{3},\frac{5\pi}{12}),\quad
(\frac{\pi}{5},\frac{\pi}{3},\frac{7\pi}{15}).$$

The Veech groups of the associated translation surfaces are
$$\triangle(9,\infty,\infty),\quad
\triangle(6,\infty,\infty),\quad
\triangle(15,\infty,\infty).$$
\end{itshape}
\end{quote}

$\T(3,4,5)$ appeared in \cite{Vch89}, $\T(3,4,5)$
and $\T(3,5,7)$ in \cite{Vrb96}; $\T(2,3,4)$ first appeared in \cite{KenSmi00} 
where Kenyon and Smillie conjectured that the three mentioned triangles 
were the only lattice examples amongst all rational acute triangles, 
and gave a proof with a bound on the common denominator of the angles (10000). 
Puchta made this bound useless \cite{Pch01}.
Hubert and Schmidt gave the precise Veech groups \cite{HbtSch01}.

\subsection{Obtuse triangles}

Ward \cite{Wrd98} defines \textit{sharp triangles} of type
$(\cfrac{\pi}{m}, \cfrac{p\,\pi}{m}, \cfrac{q\,\pi}{m})$ with $p<q$ 
and $4p\leqslant m$, and proves the following.

\begin{quote}
\begin{itshape}
  A sharp triangle with $p$ or $m$ odd has the lattice property if 
  and only if $p=1$ (in which case it is isosceles).
\end{itshape}
\end{quote}

\subsubsection{Isosceles}
The following result is due to Veech.
\begin{quote}
\begin{itshape}
  The obtuse isosceles triangle with two angles $\pi/n$ for  
  $n \geqslant 5$ has the lattice property. The Veech group of the 
  associated translation surface is $\triangle(2,n,\infty)$ 
  if $n$ is odd, and $\triangle(m,\infty,\infty)$ if $n=2m$.
\end{itshape}
\end{quote}

Hubert and Schmidt \cite{HbtSch00} show that the lattice property is 
lost in this example if we mark the points that come from the 
vertices of angle $\pi/n$ of the triangle.

\subsubsection{Scalene}

Examples by Vorobets \cite{Vrb96} and Ward \cite{Wrd98} show that 
among obtuse triangles with two angles of type $\pi/n$ and $\pi/m$, 
some have the lattice property and some do not.

In particular they studied those triangles with angles $\pi/2n$ and $\pi/n$.

\begin{quote}
\begin{itshape}
  For $n\geqslant4$, 
  the triangle $(\cfrac{\pi}{2n}, \cfrac{\pi}{n}, \cfrac{(2n-3)\pi}{2n})$ 
  has the lattice property. The Veech group of the associated surface 
  is $\triangle(3,n,\infty)$.
\end{itshape}
\end{quote}
Vorobets proved the lattice property, Ward
\cite[Theorem A]{Wrd98} identified the precise Veech group.

Vorobets shows that the triangles $\T(1,3,8)$ and $\T(2,3,7)$ do not 
have the lattice property. For $\T(1,3,8)$, it also follows from 
Ward's sharp triangles criterion.

\subsection{Final comments}
Other polygons that have been studied include rectangles, squares 
with a wall, rhombi, L-shaped polygons.

Rational billiards were the primary source of examples for Veech 
surfaces. They furnished discrete series and some isolated examples,
but in all this only gave a finite number of examples in each genus.
Other techniques recently provided infinitely many examples in genus 
two \cite{Clt, Mcm}.

\bibliographystyle{plain}

\begin{thebibliography}{wwwww99}

\bibitem[AurItz88]{AurItz88}E.~Aurell, C.~Itzykson. Rational billiards 
  and algebraic curves. \textit{J.\ Geom.\ Phys.}\ {\bf 5}:2 (1988) 191--208.
\bibitem[Brd83]{Brd83}A.~Beardon. {\em The geometry of discrete groups}. 
  GTM {\bf 91}. Springer-Verlag, 1983.\\
  NB: Corrected reprint, 1995.
\bibitem[Clt]{Clt}K.~Calta. Veech surfaces and complete 
  periodicity in genus 2.\\
  arXiv:math.DS/0205163
\bibitem[EarGrd97]{EarGrd97}C.~Earle, F.~P.~Gardiner. Teichm\"uller disks 
  and Veech's F-structures. {\em Extremal Riemann Surfaces}. Contemp.\ Math.\ 
  {\bf 201}. AMS, 1997, pp.\ 165--189.
\bibitem[FoxKsh36]{FoxKsh36}R.~H.~Fox, R.~B.~Kershner.
  Geodesics on a rational polyhedron.
  \textit{Duke Math.\ J.}\ {\bf 2} (1936) 147-150.
\bibitem[GutJdg96]{GutJdg96}E.~Gutkin, C.~Judge. The geometry and 
  arithmetic of translation surfaces with applications to polygonal 
  billiards. \textit{Math.\ Res.\ Lett.}\ {\bf 3}:3 (1996) 391--403.
\bibitem[GutJdg00]{GutJdg00}E.~Gutkin, C.~Judge. Affine mappings 
  of translation surfaces: geometry and arithmetic. \textit{Duke Math.\ J.}\ 
  {\bf 103}:2 (2000) 191--213.
\bibitem[HbtSch00]{HbtSch00}P.~Hubert, T.~A.~Schmidt. Veech groups and 
  polygonal coverings. \textit{J.\ Geom.\ Phys.}\ {\bf 35}:1 (2000) 75--91.
\bibitem[HbtSch01]{HbtSch01}P.~Hubert, T.~A.~Schmidt. Invariants of 
  translation surfaces. \textit{Ann.\ Inst.\ Fourier (Grenoble)} 
  \textbf{51}:2 (2001) 461--495.
\bibitem[KtkZml75]{KtkZml75}A.~B.~Katok, A.~N.~Zemlyakov.
  Topological transitivity of billiards in polygons.
  \textit{Math.\ Notes} {\bf 18}:2 (1975) 760-764.\\
  See also \textit{Math.\ Notes} {\bf 20}:6 (1976) 883.
\bibitem[KenSmi00]{KenSmi00}R.~Kenyon, J.~Smillie. Billiards on 
  rational-angled triangles. \textit{Comment.\ Math.\ Helv.}\ {\bf 75}:1 
  (2000) 65--108.
\bibitem[Mcm]{Mcm}C.~T.~McMullen. Billiards and Teichm\"uller curves on 
  Hilbert modular surfaces. To appear.
\bibitem[Pch01]{Pch01}J.-C.~Puchta. On triangular billiards. 
  \textit{Comment.\ Math.\ Helv.}\ \textbf{76}:3 (2001) 501--505. 
\bibitem[Vch89]{Vch89}W.~A.~Veech.  Teichm\"uller curves in moduli space, 
  Eisenstein series and an application to triangular billiards. 
  \textit{Invent.\ Math.}\ {\bf 97}:3 (1989) 553--583.\\
  {\em Erratum}: \textit{Invent.\ Math.}\ {\bf 103}:2 (1991) 447.
\bibitem[Vch92]{Vch92}W.~A.~Veech.  The billiard in a regular 
  polygon. \textit{Geom.\ and Func.\ Anal.}\ {\bf 2}:3 (1992) 341--379.
\bibitem[Vrb96]{Vrb96}Ya.~B.~Vorobets. Planar structures and billiards 
  in rational polygons: the Veech alternative. \textit{Russ.\ Math.\ Surv.}\ 
  {\bf 51}:5 (1996) 779--817.
\bibitem[Wrd98]{Wrd98}C.~C.~Ward. Calculation of Fuchsian groups 
  associated to billiards in a rational triangle. \textit{Erg.\ Th.\ 
  \& Dyn.\ Sys.}\ {\bf 18}:4 (1998) 1019--1042.
\end{thebibliography}

\medskip
\begin{small}
Samuel Leli\`evre\\
Irmar, Campus de Beaulieu, Universit\'e de Rennes 1,
35042 \textsc{Rennes}, France\\
Telephone: (+33) 223 23 58 55 \quad Fax: (+33) 223 23 67 90\\
E-mail: \texttt{samuel.lelievre@univ-rennes1.fr}\\
Web page: \texttt{http://www.maths.univ-rennes1.fr/\~{}slelievr/}
\end{small}
\end{document}